\documentclass[10pt,twocolumn,twoside]{IEEEtran}
\usepackage{bm}   % for bold
\usepackage{amsmath}
\allowdisplaybreaks[3]
\usepackage{amssymb}
\usepackage{indentfirst}  % for indent
\usepackage{graphicx}  % for pics
\usepackage{subfigure}
\usepackage{multirow} % for tabs
\usepackage{booktabs} % for tabs
\usepackage{cite}  % for citation and reference
\usepackage{threeparttable}   % for annotation in tabs
\usepackage{algorithm}% http://ctan.org/pkg/algorithm
\usepackage{algpseudocode}% http://ctan.org/pkg/algorithmicx
\usepackage{arydshln}
\usepackage{color}
\usepackage[american]{babel}
\usepackage{microtype}
\usepackage[numbers,sort&compress]{natbib}  % to include multi-citation in one []
\usepackage{mathtools}
\usepackage{epstopdf}

\def\mV{\mathcal{V}}
\def\mE{\mathcal{E}}

\def\cov{\bm{\Sigma}}
\def\covsqrt{\bm{\Sigma}^{\frac{1}{2}}}

\usepackage[font=footnotesize,justification=centering]{caption}
\def\bu{\bm{u}}
\def\bb{\bm{b}}

\def\bs{\bm{s}}
\def\bx{\bm{x}}
\def\bA{\bm{A}}

\def\bw{\bm{w}}

\begin{document}
\title{Inverse Chance Constrained Optimal Power Flow}

\author{Shenglu~Wang,
        Kairui~Feng,~\IEEEmembership{Member,~IEEE,}
        Mengqi~Xue,~\IEEEmembership{Member,~IEEE,}
        and~Yue~Song,~\IEEEmembership{Member,~IEEE}
%\vspace{-10pt}

%\thanks{This work was supported by the HKU Seed Fund for Basic Research for New Staff under Project No. 202009185007.}
% \thanks{The authors are with the Department of Control Science and Engineering, Tongji University, Shanghai 201804, China, also with the National Key Laboratory of Autonomous Intelligent Unmanned Systems, Shanghai 201210, China, and also with the Frontiers Science Center for Intelligent Autonomous Systems, Ministry of Education, Shanghai 200120, China (e-mail: ysong@tongji.edu.cn).}}

\thanks{The authors are with the National Key Laboratory of Autonomous Intelligent Unmanned Systems, Tongji University, Shanghai, China (e-mail: ysong@tongji.edu.cn).}}

\markboth{} %, ~Vol.~, No.~, ~2015}
%\markboth{network-based small-disturbance stability-Yue Song}
{Wang \MakeLowercase{\textit{et al.}}: }

\maketitle

\begin{abstract}
  The chance constrained optimal power flow (CC-OPF) essentially finds the low-cost generation dispatch scheme ensuring operational constraints are met with a specified probability, termed the security level.
  While the security level is a crucial input parameter, how it shapes the CC-OPF feasibility boundary has not been revealed.
  Changing the security level from a parameter to a decision variable,
  this letter proposes the inverse CC-OPF that seeks the highest feasible security level supported by the system.
  To efficiently solve this problem, we design a Newton-Raphson-like iteration algorithm leveraging the duality-based sensitivity analysis of an associated surrogate problem.
  Numerical experiments validate the proposed approach, revealing complex feasibility boundaries for security levels that underscore the importance of coordinating security levels across multiple chance constraints.
\end{abstract}


% \fi
\IEEEpeerreviewmaketitle

%IEEEhowto:kopka

% 1. 3.A, lambda & D_beta
% 2. 3.B, Algorithm compress
% 3. 4, Case study

\section{Introduction}\label{secintro}

Chance-constrained optimal power flow (CC-OPF) has been widely used to address the impact of uncertainties in power systems.
Over the past decade, it has evolved into a family of formulations, including Gaussian-based~\cite{bienstock2014chance}, non-Gaussian-based~\cite{wang2017chance}, and distributionally robust variants~\cite{yang2021tractable}.
Essentially, CC-OPF incorporates uncertain factors $\bw$ (e.g., wind or solar power) modeled as random variables and
seeks decision variables $\bx$ (e.g., generation dispatch) that satisfy operational security constraints with high probability,
i.e., $\Pr\{g(\bx,\bw)\leq 0\} \geq \beta$ where $\beta$ is referred to as the \textit{security level}  that serves as a key input parameter.
Therefore, the CC-OPF solutions are inherently tied to the chosen security level.

If we change the role of the security level from an input parameter to a decision variable,
a fundamental question naturally arises:
what is the maximum security level under which the CC-OPF remains feasible?
To the authors' knowledge, this question has not been systematically studied in the literature.

In this paper, we establish the inverse CC-OPF (ICC-OPF) modeling that finds the maximum security level along a given security increase direction.
To enable efficient computation of the ICC-OPF solution, we propose a Newton-Raphson-like (NR-like) iteration algorithm by leveraging the duality-based sensitivity of the optimal value of a surrogate problem.
Under the Gaussian uncertainty, we derive a convex surrogate problem and a closed-form NR-like iteration scheme which is computationally friendly.
From the numerical experiments on two IEEE test systems, the ICC-OPF provides new insights into how the maximum security levels of different chance constraints are restricting each other, which serves as a useful complement to CC-OPF.

\section{Problem Formulation}

\subsection{A Quick Review of CC-OPF}
The CC-OPF problem under a linear power flow model (e.g.,
DC power flow equation) takes the following compact form:\begin{subequations}\label{compactCCOPF}
  \begin{align}
        \min_{\bx}~&{C}(\bx)  \label{compactObj} \\
        s.t.~& \bA \bx = \bb,~\bm{E} \bx \leq \bm{f} \label{compactDeter} \\
              &\Pr\left\{\bw_{k}^T \bx \leq d_{k} \right\}\geq \beta_{k},~\forall k\in\mV_g \cup \mE \label{compactProb}
  \end{align}
\end{subequations}
where $\bx$ denotes the decision variables (e.g., active generations and participation factors of dispatchable generators);
$\mV_g$ denotes the set of generator buses and $\mE$ denotes the set of transmission lines;
$\bA$, $\bb$, $\bm{E}$, and $\bm{f}$ are the associated coefficient matrices given by the power transfer distribution factors and power flow limits;
\eqref{compactDeter} refers to deterministic constraints independent of uncertainty (e.g., power flow equation under the forecast scenario);
\eqref{compactProb} refers to chance constraints on power generations and line flows, with some entries of vector $\bw_{k}$ following certain probability distributions
and $\beta_k$ being a predefined security level to regulate the probability of constraint satisfaction.

% \vspace{-1pt}

\subsection{Formulating the ICC-OPF}

In (2) the security level $\beta_k$ is usually predefined by engineering heuristics. On the other hand,
the CC-OPF will be inevitably infeasible in case that those uncertainty factors follow long-tail distributions and $\beta_k$ is sufficiently close to 1.
Therefore, it is of interest to reconsider the CC-OPF problem in an inverse manner.
Instead of a predefined constant $\beta_k$,
we set $\beta_k$ as a variable and find the maximum security level under which problem \eqref{compactCCOPF} maintains feasible, aiming to characterize the highest security level that the system can achieve.

Thus, we propose to formulate the ICC-OPF along a certain security level increase direction as follows:\begin{subequations}\label{InverseCCOPF}
  \begin{align}
        \beta_{\max} = \max_{\bx, \beta}~&  \beta  \label{InverseObj} \\
        s.t.~  & \bA \bx = \bb,~\bm{E} \bx \leq \bm{f} \label{InverseFeasible} \\
               & \Pr\left\{\bw_{k}^T \bx \leq d_{k} \right\}\geq \beta_k  ,~\forall k\in\mV_g \cup \mE \label{InverseConstrain} \\
               & \beta_k =  \beta u_k + \beta_{0,k}, ~\forall k\in\mV_g \cup \mE \label{InverseBeta}
  \end{align}
\end{subequations}
where a common scalar variable $\beta$ scales the security level $\beta_k$ for each chance constraint along
a given unit direction vector $\bu=[u_k]$ starting at a given offset $\beta_{0,k}$.
The goal of \eqref{InverseCCOPF} is to find the largest feasible security level $\beta_{\max}$ in the direction of $\bu$, while ensuring the original CC-OPF remains feasible.

\section{Solution Method}
\subsection{The Surrogate Model with Slack Variables}
The ICC-OPF is hard to solve directly due to its highly nonlinear nature.
Alternatively, the following surrogate problem with the slack variable $\bs$ is considered:
\begin{subequations}\label{slackCCOPF}
  \begin{align}
        \|\bs\|_2^* & = \min_{\bx, \bs}~ \| \bs \|_2  \label{slackObj} \\
        s.t.~  & \eqref{InverseFeasible}, \eqref{InverseBeta}\\
               & \Pr\left\{\bw_{k}^T \bx  - s_k\leq d_{k} \right\}\geq \beta_{k} : \lambda_k, ~\forall k\in\mV_g \cup \mE \label{slackConstrain}\\
               & s_{k} \geq 0,~\forall k\in\mV_g \cup \mE  \label{slackFesible}
  \end{align}
\end{subequations}
where $ \beta $ is a given parameter and $\lambda_k \geq 0$ represents the dual variable associated with constraint \eqref{slackConstrain}.
For the surrogate problem, the optimal value $\|\bs\|_2^*$ is an implicit function of ${\beta}$.

\begin{figure}
  \centering
  \includegraphics[width=0.35\textwidth]{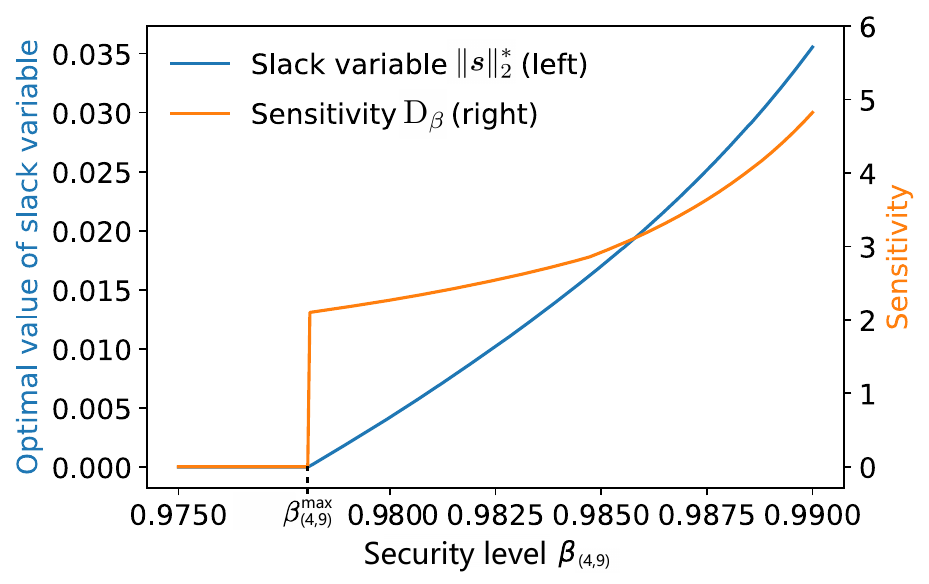}
  \caption{Illustration of the relationship between $\|\bs\|_2^*$ and $\mathrm{D}_{\beta}$ with respect to $\beta$.
  For better interpretability, the horizontal axis is linearly transformed as $\beta_{(4,9)} = \sqrt{2}\beta + 0.95$.
  The figure is plotted for the IEEE 14-bus system with parameters identical to those in Section~\ref{seccase}.
  }
  \label{fig:slackCCOPF}
\end{figure}

Problem \eqref{slackCCOPF} is always feasible since the slack variable can relax the chance constraints when necessary.
For a small $\beta_k$ with which the ICC-OPF is feasible, $\|\bs\|_2^* = 0$.
When $\beta_k$ is sufficiently large to make ICC-OPF infeasible, it follows that $\|\bs\|_2^* > 0$.
Therefore, at the value of $\beta_{\max}$, ICC-OPF is critically feasible
and the corresponding $\|\bs\|_2^*$ changes from $0$ to $0^+$, as illustrated by the blue curve in Fig. \ref{fig:slackCCOPF}.

\subsection{NR-like Iteration Scheme}
Due to the nonlinear relationship between $\|\bs\|_2^*$ and $\beta$,
a NR-like iteration scheme is designed to find $\beta_{\max}$.
The process starts with a sufficiently large
$\beta$ (such that $\beta>\beta_{\max}$ and consequently, $\|\bs\|_2^*>0$) and gradually decreases $\beta$ until the corresponding $\|\bs\|_2^*$ approaches zero.
This procedure identifies the aforementioned critical point at which $\|\bs\|_2^*=0$ but
any slight increase in $\beta$ would result in $\|\bs\|_2^*>0$.

To construct the NR-like iteration scheme,
let us begin with the first-order expansion of $\|\bs\|_2^*$ with respect to $\beta$:
\begin{equation}
  \|\bs\|_2^* = \| \bs\|_2 ^ * \big|_{\beta_0} + \mathrm{D}_{\beta}\big|_{\beta_0} (\beta - \beta_0) + \mathcal{O}((\beta - \beta_0)^2),
  \label{approximation}
\end{equation}
where $\mathrm{D}_{\beta} = \mathrm{d}\| \bs \|_2 ^ {*} \big / \mathrm{d}{\beta }$ denotes the sensitivity of the optimal objective value with respect to parameter $\beta$.
Thus, the iteration formula for computing $\beta$ towards $\|\bs\|_2^* = 0$ is:
\begin{equation}
    \beta^{(t+1)} = \beta^{(t)} - \eta_t {\| \bs \|_2 ^  {*(t)} } \big/ {\mathrm{D}_{\beta}^{(t)}},
\end{equation}
where $\beta^{(t)}$ denotes the value of $\beta$ at $t$-th iteration, and $\eta_t$ is the step size.
Meanwhile, by applying Theorem 8.2 in~\cite{conejo2006decomposition},
the sensitivity $\mathrm{D}_\beta$ can be formally expressed as the derivative of the Lagrangian function of \eqref{slackCCOPF} with respect to $\beta$,
evaluated at the optimum ${\lambda}_k^*$ and $(\bx^*,\bs^*)$. This results in:
% Meanwhile, using Theorem 8.2 in \cite{conejo2006decomposition}, we have the following formula for $\mathrm{D}_\beta$ in terms of the dual variable $\lambda_k$ in problem \eqref{slackCCOPF}:
\begin{equation}
  %\mathrm{D}_{\beta} = \frac{\partial }{\partial \beta} \mathcal{L}\left(\bx^*,\bs^*, \bm{\lambda}^*, \beta\right),
  \mathrm{D}_{\beta} = \frac{\partial }{\partial \beta}\bigg(\sum_{k\in\mV_g \cup \mE} \lambda_k^* (\beta_k  -\Pr\left\{\bw_{k}^T \bx^*  - s_k^*\leq d_{k} \right\} )\bigg) ,
  \label{D_beta}
\end{equation}
containing both explicit dependence via $\beta_k$ and implicit effects through optimal variables $\lambda_k^*, \bx^*$ and $\bs^*$ that changes with $\beta$.
%where $\mathcal{L}$ represents the Lagrangian function of the optimization problem \eqref{slackCCOPF} evaluated at the optimal solution $\bx^*,\bs^*$ and $ \bm{\lambda}^*$.
% Ignoring degenerate cases, the optimal solution $\beta_{\max}$ of the ICC-OPF problem is equivalent to satisfying $\|\bs\|_2^*= 0$ while $\mathrm{D}_{\beta} >0$.

Note that $\mathrm{D}_{\beta}$ is nonzero if and only if at least one constraint in \eqref{slackConstrain} is binding at the optimum (i.e., at least one $\lambda_k^* > 0$).
Combining this with the earlier discussion on $\beta_{\max}$ and $\|\bs\|_2^*$, we can characterize $\beta_{\max}$ as follows:
As $\beta$ keeps decreasing from a large initial value, $\|\bs\|_2^*$ reduces and eventually reaches zero at $\beta_{\max}$.
At this point, the previously binding constraints become non-binding, causing $\mathrm{D}_\beta$ to jump from a positive value to zero (see the orange curve in Fig. 1).
Therefore, $\beta_{\max}$ can be numerically identified by checking $\|\bs\|_2^* \leq \varepsilon_{s}$ and $\mathrm{D}_{\beta} > \varepsilon_D$ with prescribed tolerances $\varepsilon_s$ and $\varepsilon_D$.
Considering that the probability term in \eqref{D_beta} is still implicit, the next subsection will derive an analytical expression of $\mathrm{D}_\beta$ under the Gaussian uncertainty.

\subsection{Analytical Formulation under Gaussian Uncertainty}
Assume that the uncertain parameter $\bw_k$ follows a Gaussian distribution,
i.e., $\bw_k \sim \mathcal{N}(\bar{\bw}_{k}, \cov)$, where $\bar{\bw}_{k}$ and $\cov$ denote the mean and covariance matrix, respectively.
Moreover, in case that ${\bw}_{k}$ follows a smooth non-Gaussian distribution, it can be well approximated by a Gaussian mixture model \cite{wang2017chance},
allowing the following derivations to be conveniently extended.
Under the Gaussian uncertainty, constraint \eqref{slackConstrain} can be explicitly reformulated as \cite{bienstock2014chance}:
\begin{equation}
  \bar{\bw}_{k}^T \bx - s_k + \Phi ^{-1} (\beta_k) \|\covsqrt \bx \|_2
    \leq {d}_{k} : \lambda_k,~\forall k\in\mV_g \cup \mE  \label{formedConstrain}
\end{equation}
where $\Phi ^{-1}$ denotes the inverse cumulative distribution function (ICDF) of the standard Gaussian distribution.
Problem \eqref{slackCCOPF} becomes a second-order cone program (SOCP) when \eqref{slackConstrain} is transformed to \eqref{formedConstrain},
which can be efficiently solved by the off-the-shelf solvers and the corresponding value of the dual variable $\lambda_k$ obtained simultaneously as a byproduct.

Substituting the chance constraint expression in \eqref{formedConstrain} into \eqref{D_beta},
and applying the chain rule, we obtain
\begin{equation}
  \mathrm{D}_{\beta} =  \sum_{ k\in\mV_g \cup \mE} \frac{\partial \| \bs \|_2 ^ *}{\partial \phi_k} \frac{\mathrm{d} \phi_k}{\mathrm{d} \beta},
  \label{D_beta_chain}
\end{equation}
where $\phi_k = \Phi ^{-1} (\beta_k)$, and
\begin{equation}
  \frac{\partial \| \bs \|_2 ^ *}{\partial \phi_k} = \lambda_k^* \|\covsqrt \bx ^* \|_2, ~\forall k\in\mV_g \cup \mE.
  \label{parts_s}
\end{equation}
According to the definition of the ICDF $\Phi ^{-1}$, we have
\begin{equation}
  \frac{\mathrm{d} \phi_k}{\mathrm{d} \beta} = \frac{\mathrm{d} \phi_k}{\mathrm{d} \beta_k} u_k = \sqrt{2\pi} u_k \exp\left(\frac{\phi_k^2}{2}\right).
  \label{dphi}
\end{equation}
Substituting \eqref{parts_s}, \eqref{dphi} into \eqref{D_beta_chain} yields the computationally friendly formula for $\mathrm{D}_\beta$:
\begin{equation}
  \mathrm{D}_{\beta}
    = {\sqrt{2\pi}} \sum_{ k\in\mV_g \cup \mE}    u_k\lambda_k^* \|\covsqrt \bx^* \|_2 \exp\left(\frac{\phi_k^2}{2}\right).
  \label{sensitivity}
\end{equation}

Note that all $\lambda_k^*=0$ when $\beta<\beta_{\max}$,
which results in the zero value of $\mathrm{D}_\beta$, i.e.,  a loss of gradient.
To address this issue, a dynamic step size $\eta_t$ is designed:
starting from $\eta_t = 1$, if the updated $\beta$ leads to an interior point of the ICC-OPF ($\beta < \beta_{\max}$),
we perform a rollback and halve the step size $\eta_t$, so that each accepted $\beta$ is an exterior point of the ICC-OPF.
The detailed algorithm is presented in Algorithm 1.
\begin{algorithm}
  \caption{NR-like Solution Method for ICC-OPF}\label{al:nr}
  \begin{algorithmic}[1]
  \State Set tolerances $\varepsilon_s$, $\varepsilon_D$, iteration index $t \gets 0$, maximum number of iterations $T$, and flag $flag_{\eta} \gets \mathrm{True}$
  \State Initialize a sufficiently large  $\beta^{(t)}$ such that $\| \bs \|_2 ^ {*(t)} > 0$
  \State Solve the surrogate problem \eqref{slackCCOPF} in its SOCP form with $\beta^{(t)}$ to obtain $\| \bs \|_2 ^ {*(t)}$ , $\bm{\lambda}^ {*(t)}$
  and $\mathrm{D}_{\beta}^ {(t)}$
  \While{$t \gets t + 1< T$ and $\| \bs \|_2 ^ * > \varepsilon_s$}
    \If{$flag_{\eta}$}
      \State $\eta_t \gets 1$
    \EndIf

    \State Update the estimated value $\hat \beta \gets \beta - \eta_t {\| \bs \|_2 ^  {*} }\big/{\mathrm{D}_{\beta}}$

    \State Solve \eqref{slackCCOPF} with $\hat{\beta}$ to obtain $\|\hat{\bs}\|_2^*$, $\hat{\bm{\lambda}}^*$, and $\hat{\mathrm{D}}_{\beta}$

    \If{$| \hat {\mathrm{D}}_{\beta} | \leq \varepsilon_D$ }
      \State Set $\eta_t \gets \eta_t / 2$,
       $flag_{\eta} \gets \mathrm{False}$
    \Else
      \State $\beta \gets \hat \beta$,
      $\mathrm{D}_{\beta} \gets \hat {\mathrm{D}}_{\beta}$, $\|  {\bs} \|_2 ^{*} \gets \|  \hat{\bs} \|_2 ^ *$,
      $flag_{\eta} \gets \mathrm{True}$
    \EndIf
  \EndWhile
  \State \textbf{return} $\beta$
  \end{algorithmic}
\end{algorithm}

\section{Case Study}\label{seccase}
Take the IEEE 14-bus system and IEEE 39-bus system to verify the proposed method.
To introduce uncertainty, we incorporate renewable generators at selected buses, modeled as Gaussian-distributed injections with a covariance matrix $\boldsymbol{\Sigma}$. The matrix $\boldsymbol{\Sigma}$ is randomly generated as a positive semi-definite matrix and scaled to ensure diagonal elements are all 0.1.

For the IEEE 14-bus system, renewable generators are placed at buses 1, 3, 6, and 9. To push the system closer to feasibility boundaries, loads are set to twice their original values from MATPOWER~\cite{zimmerman2010matpower}, while the maximum flow limit $ f_{ij}^{\max} $ is set to 12.5\% of the line's susceptance $b_{ij} $.

For the IEEE 39-bus system, renewable generators are placed at buses 3, 8, 16, and 26, following the same distribution as in the 14-bus case. Loads are scaled by a factor of 1.25, and $ f_{ij}^{\max} $ is set to $ 0.1 b_{ij} $.
Additionally, for those single lines connecting generators to the network, we set $f_{ij}^{\max} = 0.2b_{ij}$.
% in the 39-bus system, some generators at the network periphery are connected by a single line, significantly impacting the OPF feasibility region. To mitigate this, the capacity $f_{ij}^{\max}$ of lines connecting buses $30,31,\cdots,38$ is doubled.

\subsection{Verification of NR-like Method}
The performance of Algorithm \ref{al:nr} can be characterized by the accuracy of the sensitivity $\mathrm{D}_{\beta}$ and optimal solution $\beta_{\max}$, which is verified by the following perturbation tests.

In the IEEE 39-bus system, set $\beta_i = 0.95, i \in \mV_g$,
generate $\beta_{ij}, (i,j) \in \mE $ using \eqref{InverseBeta} with $\beta_{0,k} = 0.95$, and set $u_{ij} = \sqrt{2}/2$ for lines $(2,3)$ and $(23,24)$ and zero for the remaining lines.
The analytical sensitivity computed by \eqref{sensitivity} is 12.5493,
which closely matches the finite-difference result 12.2301 (step size $10^{-3}$),
with a relative error of approximately 2.5\%,
validating the accuracy of formula \eqref{sensitivity}.

To validate optimality of the obtained $\beta_{\max}$, evaluate the feasibility of the original CC-OPF \eqref{compactCCOPF} at $\beta_{\max}$ and at $\beta_{\max} + 10^{-5}$.
It turns out that the former is feasible and the latter is not, verifying that $\beta_{\max}$ lies precisely at the feasibility boundary of the CC-OPF problem, and hence is the optimal solution of ICC-OPF.

\subsection{Visualized Analysis of Maximum Security Levels}

\begin{figure}[ht]
  \centering
  \subfigure[IEEE 14-bus system]{
    \includegraphics[width=0.225\textwidth,trim={8 0 7 0},clip]{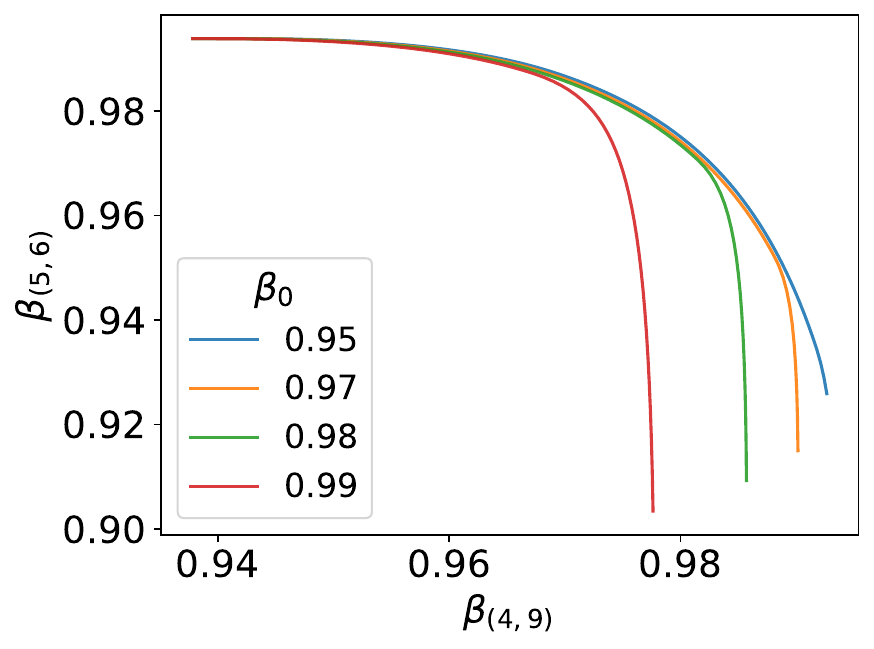}
    \label{fig:trajectory14}
  }
  \subfigure[IEEE 39-bus system]{
    \includegraphics[width=0.225\textwidth,trim={8 0 7 0},clip]{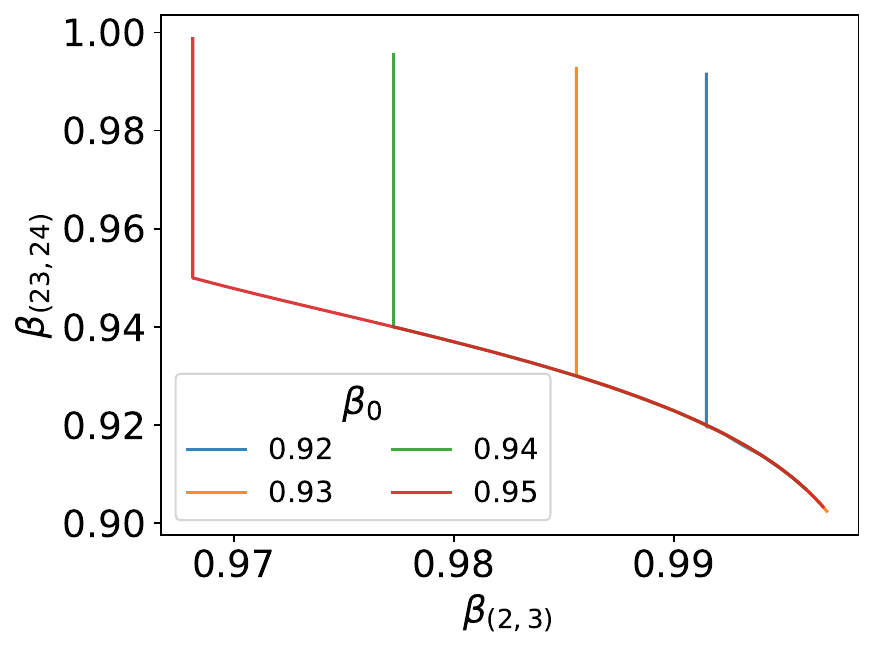}
    \label{fig:trajectory39}
  }
  \caption{Feasibility boundaries in the space of security levels.}
\end{figure}

To visualize how the maximum security levels of different lines restrain each other,
select two transmission lines in the IEEE 14-bus case and define security levels as follows:
$$
\beta_{(i,j)} =
\begin{cases}
  \beta \frac{1}{\sqrt{1+\tau^2}} + \beta_0, & (i,j) = (4,9) \\
  \beta \frac{\tau}{\sqrt{1+\tau^2}} + \beta_0, & (i,j) = (5,6) \\
  \beta_0, & \text{otherwise}
\end{cases}
$$
where $\tau$ is a tunable parameter used to generate a family of $(\beta_{(4,9)}^{\max}, \beta_{(5,6)}^{\max})$ pairs, as depicted in Fig.~\ref{fig:trajectory14}.
The lower-left region of the plot corresponds to feasible $(\beta_{(4,9)}, \beta_{(5,6)})$ pairs, while the upper-right region is infeasible.

A negative correlation is observed between $\beta_{(4,9)}^{\max}$ and $\beta_{(5,6)}^{\max}$, meaning that increasing one reduces the other.
This observation further necessitates the formulation of ICC-OPF, as the security levels of different chance constraints can significantly compete against each other and hence must be carefully set.

Furthermore, we vary the constant term $\beta_0$ to generate multiple curves. As $\beta_0$ increases (i.e., a higher security baseline for all lines), the feasible region in the lower-left compresses, which aligns with intuition.

For the IEEE 39-bus system, we select lines (2,3) and (23,24) while keeping other settings the same. The resulting curves are shown in Fig.~\ref{fig:trajectory39}.
Unlike the 14-bus case, these curves are no longer convex; in some regions, decreasing $\beta^{\max}_{(23,24)}$ does not necessarily lead to an increase in $\beta_{(2,3)}^{\max}$, revealing a more intricate non-convex relationship.
This non-convex interplay further highlights the complexity in the interaction between the maximum security levels of multiple chance constraints.

% \section{Conclusion}\label{secconclu}

% This letter proposes the inverse CC-OPF problem, which investigates the maximum security levels the system can sustain while preserving CC-OPF feasibility.
% A sensitivity-based Newton-Raphson-like algorithm is designed for efficient computation.
% Numerical results on IEEE 14-bus and 39-bus systems demonstrate the impact of transmission line security level selection on the feasibility of the overall system.

\ifCLASSOPTIONcaptionsoff
  \newpage
\fi

{\footnotesize
\bibliographystyle{IEEEtran}
\bibliography{IEEEabrv,inverse-ccopf}

\begin{thebibliography}{1}
\providecommand{\url}[1]{#1}
\csname url@samestyle\endcsname
\providecommand{\newblock}{\relax}
\providecommand{\bibinfo}[2]{#2}
\providecommand{\BIBentrySTDinterwordspacing}{\spaceskip=0pt\relax}
\providecommand{\BIBentryALTinterwordstretchfactor}{4}
\providecommand{\BIBentryALTinterwordspacing}{\spaceskip=\fontdimen2\font plus
\BIBentryALTinterwordstretchfactor\fontdimen3\font minus
  \fontdimen4\font\relax}
\providecommand{\BIBforeignlanguage}[2]{{%
\expandafter\ifx\csname l@#1\endcsname\relax
\typeout{** WARNING: IEEEtran.bst: No hyphenation pattern has been}%
\typeout{** loaded for the language `#1'. Using the pattern for}%
\typeout{** the default language instead.}%
\else
\language=\csname l@#1\endcsname
\fi
#2}}
\providecommand{\BIBdecl}{\relax}
\BIBdecl

\bibitem{bienstock2014chance}
D.~Bienstock, M.~Chertkov, and S.~Harnett, ``Chance-constrained optimal power
  flow: Risk-aware network control under uncertainty,'' \emph{{SIAM} Review},
  vol.~56, no.~3, pp. 461--495, 2014.

\bibitem{wang2017chance}
Z.~Wang, C.~Shen, F.~Liu, X.~Wu, C.-C. Liu, and F.~Gao, ``Chance-constrained
  economic dispatch with non-gaussian correlated wind power uncertainty,''
  \emph{IEEE Transactions on Power Systems}, vol.~32, no.~6, pp. 4880--4893,
  2017.

\bibitem{yang2021tractable}
L.~Yang, Y.~Xu, H.~Sun, and W.~Wu, ``Tractable convex approximations for
  distributionally robust joint chance-constrained optimal power flow under
  uncertainty,'' \emph{IEEE Transactions on Power Systems}, vol.~37, no.~3, pp.
  1927--1941, 2021.

\bibitem{conejo2006decomposition}
A.~J. Conejo, E.~Castillo, R.~M\'inguez, and R.~Garc\'ia-Bertrand,
  \emph{Decomposition Techniques in Mathematical Programming: Engineering and
  Science Applications}.\hskip 1em plus 0.5em minus 0.4em\relax Berlin,
  Heidelberg: Springer, 2006.

\bibitem{zimmerman2010matpower}
R.~D. Zimmerman, C.~E. Murillo-S{\'a}nchez, and R.~J. Thomas, ``Matpower:
  Steady-state operations, planning, and analysis tools for power systems
  research and education,'' \emph{IEEE Transactions on Power Systems}, vol.~26,
  no.~1, pp. 12--19, 2010.

\end{thebibliography}
}

% You can push biographier them. The appropriate

%\vfill


\end{document}